\documentclass{amsart}
\pdfoutput=1
\usepackage{amssymb}
\usepackage{amsmath}
\usepackage{graphicx}
\usepackage[ruled,vlined]{algorithm2e}

\newtheorem{theorem}{Theorem}
\newtheorem{definition}{Definition}
\newtheorem{lemma}[theorem]{Lemma}
\newtheorem{remark}[theorem]{Remark}
\newtheorem{corollary}[theorem]{Corollary}

\newcommand{\comment}[1]{\mbox{}}
\newcommand{\graph}[2]{(\mathcal #1, \mathcal #2)}
\newcommand{\supp}[1]{\mathrm{supp}\left(#1\right)}
\newcommand{\conefam}{\mathcal M=\left(Q,\{(M_i,C_i)\}_{i=1}^N\right)}

\begin{document}
\author{Daniel Wilczak}
\title[Hyperbolic attractor in the Kuznetsov system]{
Uniformly hyperbolic attractor of the Smale-Williams type for a Poincar\'e map in the Kuznetsov system}

\begin{abstract}
We propose a general algorithm for computer assisted verification of uniform hyperbolicity for maps which exhibit a robust attractor.

The method has been successfully applied to a Poincar\'e map for a system of coupled non-autonomous van der Pol oscillators. The model equation has been proposed by Kuznetsov \cite{K} and the  attractor seems to be of the Smale-Williams type.
\end{abstract}
\keywords{Uniform hyperbolicity, cone condition, computer
assisted proof}
\address{
    Jagiellonian University, Institute of Computer Science\\
    {\L}ojasiewicza 6, 30-348 Krak\'ow, Poland
}
\email{wilczak@ii.uj.edu.pl}
\urladdr{http://www.ii.uj.edu.pl/~wilczak}
\date{\today}
\dedicatory{To the memory of my mathematics teacher,\\ Maria R\'og}
\subjclass[2010]{37D05, 
34D45, 
37D45
}
\maketitle

\section{Introduction.}
Hyperbolic systems of dissipative type, contracting the phase space volume, manifest robust attractors. That means that the dynamics does not change qualitatively when the parameters of the system vary in some range.

There are several models of maps which produce hyperbolic nontrivial attractor
-- like the Plykin attractor \cite{P} or the Smale solenoid \cite{KH}. But for long time there were no examples of continuous systems that apparently have hyperbolic strange attractors.

\begin{figure}[htbp]
    \centerline{
        \includegraphics[height=2in]{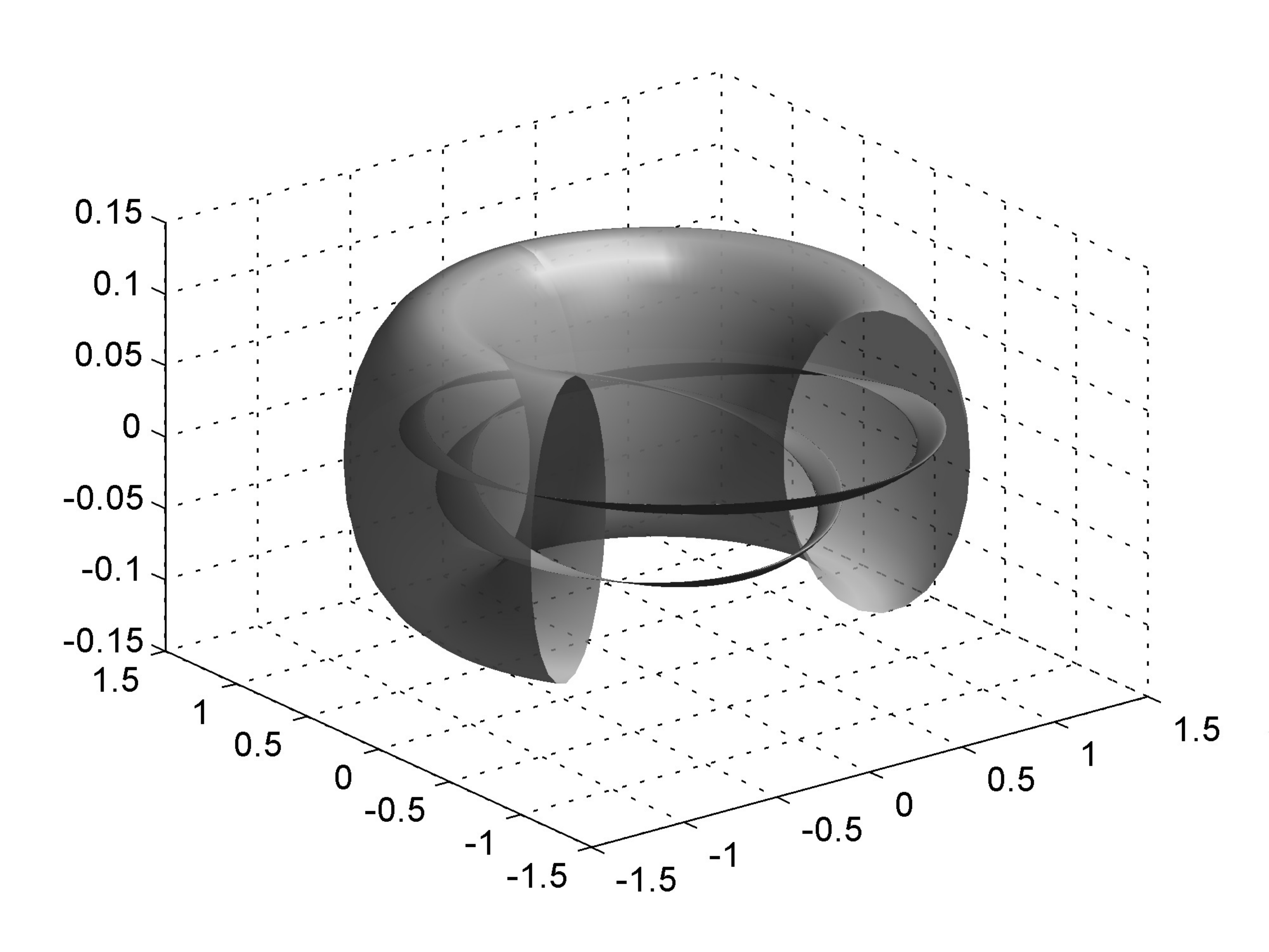}
        \includegraphics[height=2in]{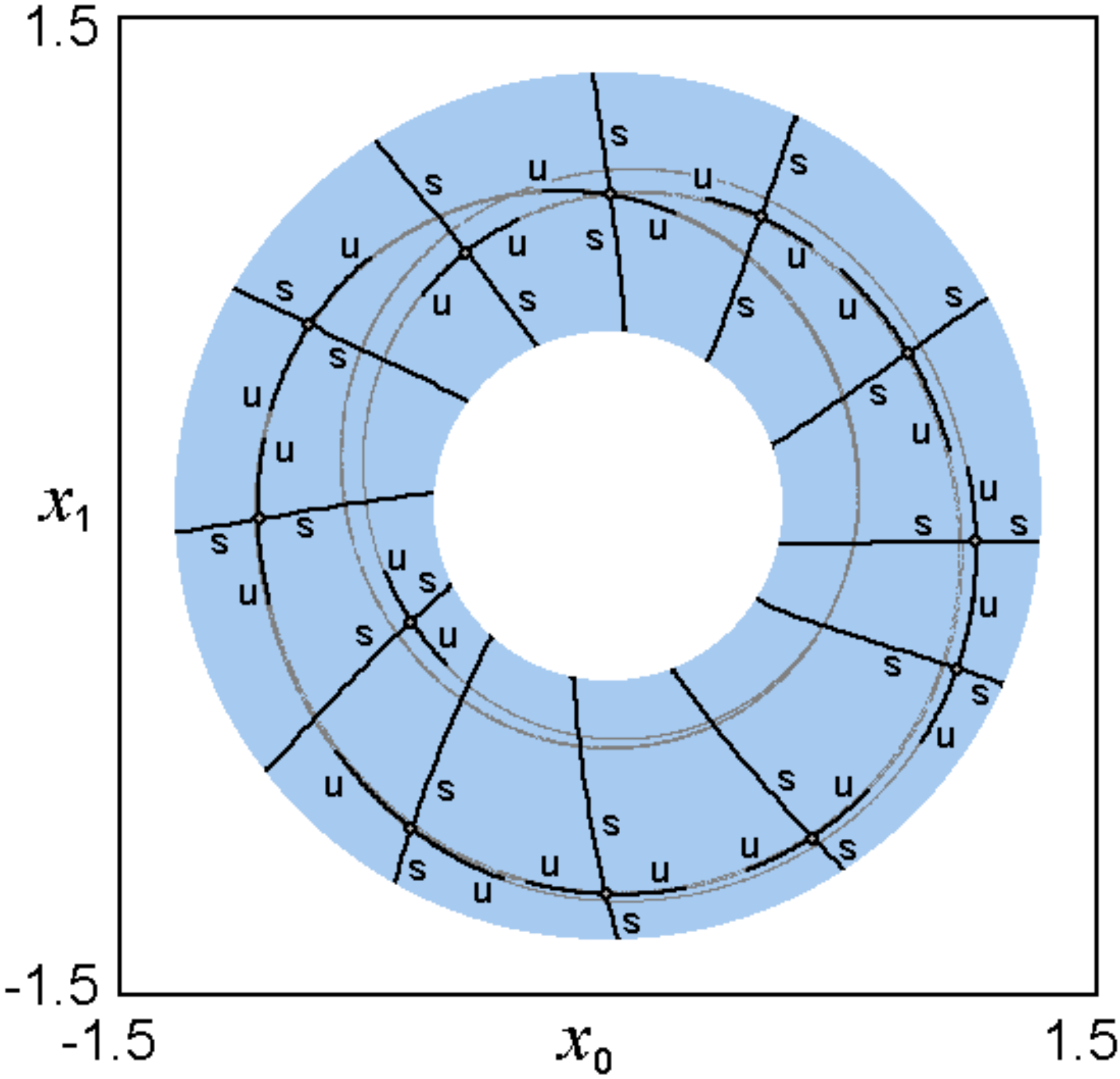}
    }
    \caption{Left: an absorbing domain projected onto three coordinates after a linear change of coordinates.
             Right: local unstable and stable foliations. It is a
             numerical indication that they are not tangent over the absorbing domain.
            These figures are \cite[Figures 1 and 7]{KS} -- used here with the permission of prof. Kuznetsov.\label{fig:attractor}}
\end{figure}

Recently, Kuznetsov \cite{K} proposed a continuous model that exhibits an uniformly hyperbolic attractor of the Smale-Williams type in the Poincar\'e map. This is a non-autonomous system in $\mathbb R^4$ and it is constructed on a basis of two coupled van
der Pol oscillators:
\begin{equation}\label{eq:vpo}
\begin{cases}
    \dot x &= \omega_0u,\\
    \dot u &=-\omega_0x + \left(A\cos(2\pi t/T ) -x^2 \right)u +(\varepsilon /\omega_0)y
                \cos(\omega_0 t),\\
    \dot y &= 2\omega_0v,\\
    \dot v &=-2\omega_0y + \left( -A\cos(2\pi t/T ) -y^2 \right) v
            + (\varepsilon /2\omega_0)x^2.
\end{cases}
\end{equation}
Let $P$ denote the Poincar\'e map of the above system defined naturally as the shift along the trajectories over the period of the vector field
\begin{equation}\label{eq:poincare}
 P(x,u,y,v) = \left(x(T),u(T),y(T),v(T)\right).
\end{equation}
Kuznetsov and Sataev \cite{KS} gave a deep numerical study of this system and observed that for some range of parameter values there is an absorbing domain for $P$ of the toroidal shape as presented in
Figure \ref{fig:attractor}, left panel. This domain is a product of
a 3D ball and a circle.

The observed attractor is of the Smale-Williams type. Numerical studies \cite{K,KS,KSe} gave a strong numerical evidence of the existence of
uniformly hyperbolic attractor for some range of parameter values
around
\begin{equation}\label{eq:parameters}
    \omega_0 = 2\pi,\quad A=5, \quad T=6, \quad \varepsilon = 0.5
\end{equation}
-- see Figure~\ref{fig:attractor} right panel, where the unstable and stable foliations are presented.

In this paper we propose a method for computer assisted verification that a map possesses an uniformly hyperbolic attractor. As a test case we apply the method to the model map proposed by Smale \cite{KH}. This is a map defined on the two-dimensional solid torus and it is analytically proved to be uniformly hyperbolic.

\begin{figure}[htbp]
 \centerline{
    \includegraphics[width=0.5\textwidth]{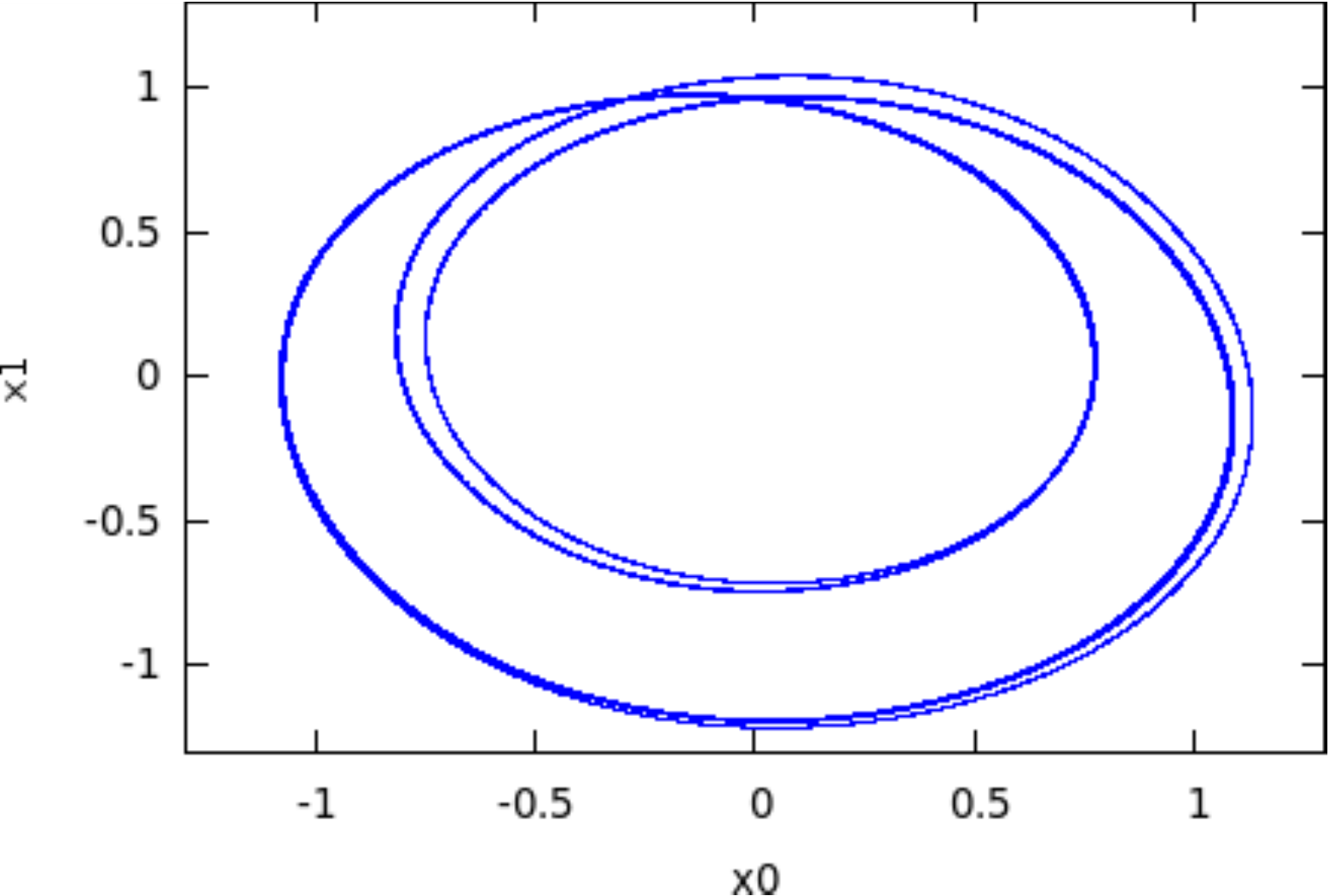}
    \includegraphics[width=0.5\textwidth]{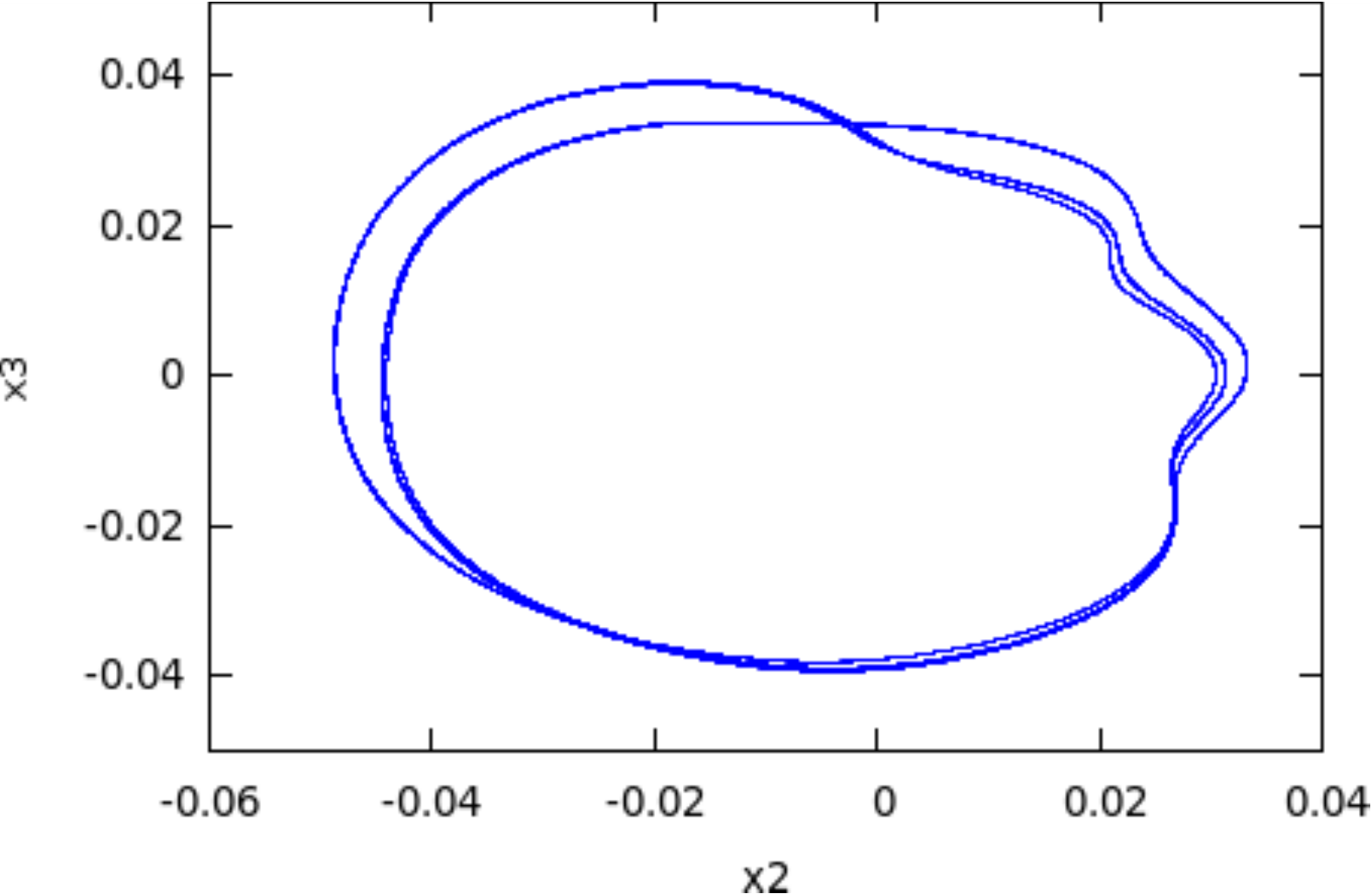}
  }
\caption{Projection of the enclosure of the attractor for $\tilde P$ computed with the resolutions $k=14$.\label{fig:sw}}
\end{figure}

Although the method is general, the main motivation for us to undertake this study was to prove Kuznetsov's conjecture about the hyperbolicity of the system (\ref{eq:vpo}). The following theorem is the main result of this paper.
\begin{theorem}\label{thm:main}
Consider the system (\ref{eq:vpo}) with the parameter values (\ref{eq:parameters}). Let $P$ denote the Poincar\'e map for this system as defined in (\ref{eq:poincare}).
There exists a compact, connected and explicitly given set $\mathcal B$ such that
\begin{enumerate}
\item $\mathcal B$ is positive invariant with respect to $P$, i.e. $P(\mathcal B)\subset \mathcal B$,
\item $P$ is uniformly hyperbolic on the maximal invariant set $\mathcal A = \bigcap_{i>0}P^i(\mathcal
B)$ with one positive and three negative Lyapunov exponents.
\end{enumerate}
\end{theorem}
After a linear change of coordinates the set $\mathcal B$ is a union of $7\,970\,392$ four-dimensional explicitly given cubes. The set $\mathcal B$ is a very narrow enclosure of $\mathcal A$ as presented in Figure~\ref{fig:sw}. Note, the coordinates $(x_0,x_1,x_2,x_3)$ are related to $(x,u,y,v)$ by a linear change of coordinates.

The above theorem implies that the set $\mathcal A = \bigcap_{i>0}P^i(\mathcal B)$ is compact and connected. The following theorem implies that it is a non-trivial continuum.
\begin{theorem}\label{thm:nontrivAttractor}
The set $\mathcal A$ contains at least one fixed point and one period-two orbit for $P$.\end{theorem}

In the literature there are already algorithms for computer assisted verification of uniform hyperbolicity and enclosure of attractors. The very pioneer and famous work has been done by Warwick Tucker \cite{T}. He proved that the Lorenz system for classical parameter values satisfies the geometric model by Guckenheimer and Holmes \cite{GH}. Moreover, he gave the proof of the existence of the SRB measure supported on the attractor.

Recently, Hruska \cite{Hr1,Hr2} proposed a method for computer assisted verification of uniform hyperbolicity. She successfully applied the method to the complex H\'enon map. The method proposed by Hruska uses the notion of box-hyperbolicity. This method reduces the verification of the hyperbolicity of an invariant set to the verification if some quadratic forms are positive definite. The box-hyperbolicity requires some conditions for the derivative of a map under consideration and its inverse. This is not a limitation in the case of the H\'enon map but makes the algorithm difficult to apply for ODE's. To compute the derivative of the inverse of a Poincar\'e map we can either invert an interval matrix that is an enclosure of the derivatives of the Poincar\'e map or integrate the ODE backwards together with the variational equations. In the first case the derivatives of Poincar\'e maps are often non-invertible as interval matrices because of unavoidable over-estimations when integrating the variational equations. Backward integration does not help in many cases. After changing the time $t\to -t$ in the equation the system often becomes stiff and again computationally difficult. In particular, this happens when the invariant set is an attractor with strong dissipation.

Later, Arai \cite{A} proposed a method for verification that a chain recurrent set for a map is uniformly hyperbolic. He applied the method to the H\'enon map \cite{H} and verified that it is hyperbolic on the non-wandering set for a wide range of parameter values. This method, however, requires huge memory to work and it is computationally expensive by its very construction. Therefore, there is a little hope to apply it successfully in the high dimensional space and for nontrivial ODE's.

Similar results were obtained by Mazur, Tabor and Ko\'scielniak \cite{MTK} and by Mazur and Tabor \cite{MT}. The authors introduced a notion of semi-hyperbolicity and applied the method to the real H\'enon map.

Our method is similar in the spirit to that proposed by Hruska. We use a notion of strong hyperbolicity proposed in \cite{KWZ,Z} as a theoretical tool for verifying uniform hyperbolicity. The main difference of the method proposed here is that it does not involve conditions for the inverse of the map under consideration.

We believe that a successful application of our algorithms to a four dimensional Poincar\'e map is a serious test for the method as well as the implementation and it proves their applicability. We would like to mention here that verification of uniform hyperbolicity for explicit given maps is significantly easier than for the continuous systems. It is clear that in a numerical approach to hyperbolic dynamics one needs to enclose both values and derivatives of the map under consideration. This is an easy task when the map is given explicitly, while for Poincar\'e maps one needs to integrate the system with its variational equations.

For rigorous integration of the system and its partial derivatives we used $C^0$ and $C^1$ solvers from the CAPD library \cite{CAPD} that the author is one of the main developers.

The paper is organized as follows. In Section~\ref{sec:theory} we give the notion of strong hyperbolicity and the theoretical background on how to verify uniform hyperbolicity using this tool. In Section~\ref{sec:algorithms} we present algorithms used to define stable and unstable cones over the attractor. In Section~\ref{sec:applications} we give an application of the proposed method to the Smale map.  We also present proofs of Theorem~\ref{thm:main} and Theorem~\ref{thm:nontrivAttractor}.

\section{Theoretical background.}\label{sec:theory}
Let $f\colon\mathbb R^n\to\mathbb R^n$ be a diffeomorphism and let $M\subset \mathbb R^n$ be a compact invariant set for $f$. We denote by $TM$ the restriction of the tangent bundle $T\mathbb R^n$ to $M$.
\begin{definition}
 $f$ is uniformly hyperbolic on $M$ if $TM$ splits into a direct sum $TM=E^u\oplus E^s$ of two $Tf$-invariant sub-bundles and there are constants $c>0$ and $0<\lambda<1$ such that
\begin{equation*}
 \|Df^n|_{E^s} v\|<c\lambda^n\|v\|\quad\text{and}\quad \|Df^{-n}|_{E^u} v\|<c\lambda^n \|v\|
\end{equation*}
hold for $n\geq0$.
\end{definition}
We will recall and reformulate the definition of strong hyperbolicity introduced in \cite{KWZ} and later extended in \cite{Z}.

Let $M=\bigcup_{i=1}^N M_i$ where $M_i$ are compact sets having pairwise disjoint interiors. In our algorithms these sets will be boxes in some coordinate systems. Let us fix nonegative integers $u,s$ such that $u+s=n$. Assume that at each set $M_i$ we have fixed a linear coordinate system $C_i$. We define a quadratic form
\begin{equation}\label{eq:Qform}
 Q(x,y) = \|x\|^2-\|y\|^2, \quad x\in\mathbb R^u, \ y\in\mathbb R^s.
\end{equation}
For $i=1,\ldots,N$ we define positive and negative cones by
\begin{eqnarray*}
 Q^+_i =\left\{u\in\mathbb R^n : Q(C_iu)>0\right\},\\
 Q^-_i =\left\{s\in\mathbb R^n : Q(C_is)<0\right\}.
\end{eqnarray*}
We will denote this family by $\conefam$ and call \emph{cubical set with cones}.

Let  $\conefam$ be a cubical set with cones and put $M=\bigcup_{i=1}^N{M_i}$. For $i,j=1\ldots,N$ we put $f_{ij}=C_j f C^{-1}_i$.

\begin{definition}
$f$ is strongly hyperbolic on $\conefam$ if for $z\in M_i$ and $j=1,\ldots,N$ such that $f(M_i)\cap M_j\neq\emptyset$ the matrix
\begin{equation*}
 [C_jDf(z)C_i^{-1}]^T Q [C_jDf(z)C_i^{-1}]-Q
\end{equation*}
is positive definite.
\end{definition}

By $\mathrm{Inv}(f,M)$ we will denote the maximal invariant set for $f$ in $M$, i.e.
\begin{equation*}
 \mathrm{Inv}(f,M) = \left\{x\in M:f^n(x)\in M, f^{-n}(x)\in M,\ \text{for }n\geq0\right\}.
\end{equation*}

\begin{theorem}\label{thm:verifyConesByQuadForm}
If $f$ is strongly hyperbolic on $\conefam$ then $f$ is uniformly hyperbolic on $\mathcal H=\mathrm{Inv}(f,M)$.
\end{theorem}

The proof of the above theorem consists of several lemmas.
\begin{lemma}\label{lem:mapCones}
 Assume $f$ is strongly hyperbolic on $\conefam$.
\begin{enumerate}
\item If $f(M_i)\cap M_j\neq\emptyset$ then for $z\in M_i$ we have $Df(z)Q_i^+\subset Q^+_j$.
\item For $z\in f(M_i)\cap M_j$ we have $(Df(z))^{-1}Q_j^-\subset Q^-_i$.
\end{enumerate}
\end{lemma}
\begin{proof}
We will prove the first assertion. Since $u\in Q^{+}_i$ we have $Q(C_iu)>0$. Fix $j$ such that $f(M_i)\cap M_j\neq\emptyset$. Put $\bar u = C_i u$ and $\bar z = C_i z$. Since $C_i$ and $C_j$ are linear and $f$ is strongly hyperbolic on $\mathcal M$ we have
\begin{equation*}
 Q(C_j Df(z)u) = Q\left(Df_{ij}(\bar z)\bar u\right) > Q(\bar u)  = Q(C_iu)>0.
\end{equation*}
Thus, by definition $Df(z)u\in Q^+_j$.

We will now prove the second assertion. Fix $z\in f(M_i)\cap M_j$ and put $\bar z = f^{-1}(z)\in M_i$. Assume that $(Df(z))^{-1}Q_j^-\not\subset Q^-_i$, hence there is $s\in Q_j^-$ such that $\bar s := (Df(z))^{-1}s\notin Q^-_i$. Reasoning as in the proof of the first assertion we conclude that
\begin{equation*}
 Q(C_j Df(\bar z)\bar s) > Q(C_i\bar s) \geq 0.
\end{equation*}
Hence, $s=Df(\bar z)\bar s\in Q^+_j$ which contradicts the choice of $s$.
\end{proof}

\begin{lemma}\label{lem:openSetOfQuadForms}
Assume $f$ is strongly hyperbolic on $\conefam$. There is $\varepsilon>0$ such that for $\lambda\in (1-\varepsilon,1+\varepsilon)$, $i,j=1,\ldots,N$ such that $f(M_i)\cap M_j\neq\emptyset$ and $z\in M_i$ we have
\begin{equation}\label{eq:epsDef}
 [C_jDf(z)C_i^{-1}]^T Q [C_jDf(z)C_i^{-1}]-\lambda Q.
\end{equation}
\end{lemma}
\begin{proof}
 The set of positive definite matrices is open. Thus for fixed $i,j$ and $z\in M_i$ we can find a neighborhood $V_{i,j,z}\subset\mathbb R^n\times \mathbb R$ of $(z,1)$ such that for $(\bar z, \lambda)\in V_{i,j,z}$ the matrix
\begin{equation*}
 [C_jDf(\bar z)C_i^{-1}]^T Q [C_jDf(\bar z)C_i^{-1}]-\lambda Q
\end{equation*}
is positive definite. Then the assertion follows from the compactness of $M$.
\end{proof}

Before we state the next lemmas let us define some constants. Put
\begin{eqnarray}
D_1(\mathcal M) &=& \max\{\|C_i\| : i=1,\ldots,N\}, \label{eq:D1Def}\\
D_2(\mathcal M) &=& \max\{\|C_i^{-1}\| : i=1,\ldots,N\}. \label{eq:D2Def}
\end{eqnarray}

\begin{lemma}\label{lem:normU_lowerBound}
Let $\conefam$. There is a constant $R>0$ such that for $u\in\mathbb R^n$ and $i=1,\ldots, N$ we have $\|u\|^2\geq R|Q(C_iu)|$.
\end{lemma}
\begin{proof}
 We have
\begin{multline*}
 \|u\|^2 = \|C_i^{-1} C_i u\|^2 \geq \|C_i\|^{-2}\|C_iu\|^2 \geq \\
 \|C_i\|^{-2}\|\pi_x C_i u\|^2\geq
 \|C_i\|^{-2}Q(C_iu) \geq
D_1(\mathcal M)^{-2}Q(C_iu).
\end{multline*}
Similarly
\begin{equation*}
 \|u\|^2 \geq
\|C_i\|^{-2}\|\pi_y C_i u\|^2 \geq -D_1(\mathcal M)^{-2}Q(C_iu).
\end{equation*}
Hence, the assertion follows with $R=D_1(\mathcal M)^{-2}$.
\end{proof}

\begin{lemma}\label{lem:expansionPositiveCones}
Assume $f$ is strongly hyperbolic on $\conefam$. There are constants $\lambda>1$, $c>0$ such that for all $i=1,\ldots,N$, $z\in M_i\cap \mathrm{Inv}(f,M)$, $u\in Q_i^+$ and $k>0$ we have $\|Df^k(z)u\|\geq c\lambda^k\|u\|$.
\end{lemma}
\begin{proof}
From Lemma~\ref{lem:openSetOfQuadForms} there is a constant $\bar\lambda>1$ such that the matrices
\begin{equation}\label{eq:modifiedSequence}
V_{i,j,z}=[C_j Df(z) C_i^{-1}]^T Q[C_j Df(z) C_i^{-1}] - \bar\lambda Q
\end{equation}
are positive definite
for $i,j=1,\ldots,N$ such that $f(M_i)\cap M_j\neq \emptyset$ and $z\in M_i$. Since $M$ is compact there is $L>0$ such that
\begin{equation}\label{eq:Vnorm1}
 V_{i,j,z}(u)\geq L\|u\|^2,\quad u\in\mathbb R^n
\end{equation}
for $i,j=1,\ldots,N$ such that $f(M_i)\cap M_j\neq \emptyset$ and $z\in M_i$.

Fix  $z=\mathrm{Inv}(f,M)$. Put $z_k = f^k(z)$ and $A_k = Df(z_k)$ for $k\geq 0$. Let us fix a sequence $\{i_k\}_{k\geq 0}$ such that $z_k\in M_{i_k}$. Note, this sequence might be not unique. Put $\tilde A_k = C_{i_{k+1}} A_k C_{i_k}^{-1}$, $k\geq 0$. From (\ref{eq:modifiedSequence}) for $u\neq 0$ and $k>1$ we have
\begin{multline}\label{eq:coneSequence}
 Q(C_{i_k}Df^k(z_0)u) = Q(\tilde A_{k-1}\cdots \tilde A_0 C_{i_0}u) > \\ \bar\lambda Q(\tilde A_{k-2}\cdots \tilde A_0 C_{i_0}u)>\cdots > \bar\lambda^{k-1} Q(\tilde A_0C_{i_0}u).
\end{multline}

From Lemma~\ref{lem:normU_lowerBound} we get that there is a constant $R>0$ depending on $\mathcal M$ only, such that for $u\in\mathbb R^n$ holds
\begin{equation}\label{eq:Dfu_lowerBound}
 \|Df^k(z_0)u\|^2 \geq R |Q(C_{i_k}Df^k(z_0)u)|.
\end{equation}
From (\ref{eq:modifiedSequence}--\ref{eq:Vnorm1}) for $u\in Q_{i_0}^+$  we have
\begin{multline}\label{eq:u_upperBound}
  Q(\tilde A_0 C_{i_0}u)\geq Q(C_{i_0}u) + L\|C_{i_0}u\|^2>L\|C_{i_0}u\|^2\geq \\ L \|C_{i_0}^{-1}\|^{-2}\|u\|^2\geq LD_2(\mathcal M)^{-2}\|u\|^2,
\end{multline}
where $D_2(\mathcal M)$ was defined in (\ref{eq:D2Def}). Combining (\ref{eq:coneSequence}--\ref{eq:u_upperBound}) we obtain that for $u\in Q_{i_0}^+$ and  $k>1$ holds
\begin{equation*}
 \|Df^k(z_0)u\|\geq c\lambda^k\|u\|
\end{equation*}
with $\lambda = \bar\lambda^{1/2}$ and $c=(R L)^{1/2} \left(\lambda D_2(\mathcal M)\right)^{-1}$. Clearly, adjusting the constant $c$ we can obtain required inequality for $k>0$.
\end{proof}

There remains for us to show backward expansion in the negative cones.
\begin{lemma}\label{lem:expansionNegativeCones}
Assume $f$ is strongly hyperbolic on $\conefam$. There are constants $\lambda>1$, $c>0$ such that for $i=1,\ldots,N$, $z\in \mathrm{Inv}(f,M)\cap M_i$, $s\in Q_i^-$ and $k>0$ we have $\|Df^{-k}(z)s\|\geq c\lambda^k\|s\|$.
\end{lemma}
\begin{proof}
From Lemma~\ref{lem:openSetOfQuadForms} there is $\bar\lambda\in(0,1)$ such that the matrices
\begin{equation}\label{eq:modifiedSequence2}
V_{i,j,z}=[C_j Df(z) C_i^{-1}]^T Q[C_j Df(z) C_i^{-1}] - \bar\lambda Q
\end{equation}
are positive definite for $i,j=1,\ldots,N$ such that $f(M_i)\cap M_j\neq \emptyset$ and $z\in M_i$. Since $M$ is compact there is $L>0$ such that
\begin{equation}\label{eq:Vnorm2}
 V_{i,j,z}(s)\geq L\|s\|^2,\quad s\in\mathbb R^n
\end{equation}
for $i,j=1,\ldots,N$ such that $f(M_i)\cap M_j\neq \emptyset$ and $z\in M_i$.

Let us fix $z\in \mathrm{Inv}(f,M)$. Choose a sequence $\{i_k\}_{k\geq 0}$ such that $z_k=f^{-k}(z)\in M_{i_k}$ for $k\geq 0$. For $k>0$ put $A_k=Df(z_k)$, $\tilde A_k = C_{i_{k-1}}A_kC_{i_k}^{-1}$ and $s_k=Df^{-k}(z_0)s$. From (\ref{eq:modifiedSequence2}) for $k>0$ and $ s\neq 0$ we obtain
\begin{equation*}
 Q(C_{i_1}Df^{k-1}(z_k)s_k) = Q(\tilde A_2\cdots \tilde A_kC_{i_k}s_k) >  \bar\lambda^{k-1} Q(C_{i_k}s_k),
\end{equation*}
 that after substituting $s_k$ and $z_k$ becomes
\begin{equation}\label{eq:invSequence}
 - Q(C_{i_k}Df^{-k}(z_0)s) > - \bar\lambda^{-k+1}Q(C_{i_1}s_1).
\end{equation}
From Lemma~\ref{lem:normU_lowerBound} there is a constant $R>0$ depending on $\mathcal M$ only, such that for $s\in\mathbb R^n$ holds
\begin{equation}\label{eq:minusQUpperBound}
 \|Df^{-k}(z_0)s\|^2\geq -R Q(C_{i_k}Df^{-k}(z_0)s).
\end{equation}
From (\ref{eq:modifiedSequence2}--\ref{eq:Vnorm2}) we have
\begin{equation*}
 Q(C_{i_0}s) = Q(C_{i_1}Df(z_1)s_1) = Q(\tilde A_1 C_{i_1}s_1) \geq  Q(C_{i_1}s_1) + L\|C_{i_1}s_1\|^2.
\end{equation*}
Hence, for $s\in Q_{i_0}^-$ holds
\begin{equation}\label{eq:minusQLowerBound}
 -Q(C_{i_1}s_1) \geq L\|C_{i_0}s\|^2 - Q(C_{i_0}s) >L\|C_{i_0}s\|^2\geq LD_2(\mathcal M)^{-2}\|s\|^2.
\end{equation}
Combining (\ref{eq:invSequence}--\ref{eq:minusQLowerBound}) we obtain that for $k>1$ and $s\in Q_{i_0}^-$ holds
\begin{equation*}
 \|Df^{-k}(z_0)s\|\geq c \lambda^k \|s\|,
\end{equation*}
with $\lambda = \bar\lambda^{-1/2}>1$ and $c=(LR)^{1/2}\left(\lambda D_2(\mathcal M)\right)^{-1}$. Adjusting the constant $c$ if necessary we can obtain the same inequality for $k>0$.
\end{proof}

\begin{proof}[Proof of Theorem~\ref{thm:verifyConesByQuadForm}]
Let us associate to each $z\in\mathrm{Inv}(f,M)$ an index $i_z$ such that $z\in M_{i_z}$. This index might be not unique but for each $z$ we can fix one of them. For simplicity we will write $M_z,C_z,Q_z^\pm$ instead of $M_{i_z}, C_{i_z}, Q_{i_z}^\pm$.

From Lemma~\ref{lem:mapCones}, Lemma~\ref{lem:expansionPositiveCones} and Lemma~\ref{lem:expansionNegativeCones} there are constants $c>0$ and $\lambda>1$ such that for $z\in\mathrm{Inv}(f,M)$ we have
\begin{eqnarray*}
 Df(z)Q_z^+&\subset& Q_{f(z)}^+,\\
 Df(z)^{-1}Q_z^-&\subset& Q_{f^{-1}(z)}^-,\\
 \|Df^k(z)u\|&\geq& c\lambda^k\|u\|, \quad\text{for } u\in Q_z^+,\\
 \|Df^{-k}(z)s\|&\geq& c\lambda^k\|s\|, \quad\text{for } s\in Q_z^-.
\end{eqnarray*}
Note, the cones $Q_z^+$ and $Q_z^-$ are disjoint and the sum of them spans $\mathbb R^n$. Hence, the assertion follows from \cite[Corollary 6.4.8]{KH}.
\end{proof}

\section{Algorithms.}\label{sec:algorithms}
\subsection{Graph representations of maps.}

Let $X$ be a topological space and let $N\subset X$. For a family of subsets $\mathcal U\in 2^{X}$ by $\supp{\mathcal U}$ we will denote its geometrical representation $\supp{\mathcal U}=\bigcup_{U\in\mathcal U} U$.

\begin{definition}
We will say that $\mathcal U\in 2^X$ is a cover of $N$ if $\mathcal U$ is a finite set and $N\subset\supp{\mathcal U}$.
\end{definition}

\begin{definition}
 Let $\mathcal U, \mathcal V$ be a finite sets. A multivalued function $f\colon \mathcal U\multimap\mathcal V$ is called a combinatorial map.
\end{definition}

Covers of compact sets as well as combinatorial maps appear naturally in the computer assisted proofs in dynamical systems.
We usually cannot compute the range of a map $f:X\to Y$ over a domain  $U\subset X$. Instead, one can fix covers $\mathcal X, \mathcal Y$ of $X$ and $Y$, respectively, and then compute combinatorial representation of the map in the covers $\mathcal X, \mathcal Y$.
This is usually realized in the following steps.
\begin{itemize}
\item For $U\subset X$ compute a minimal cover $$\mathcal U = \bigcap \left\{ \mathcal V\subset \mathcal X : U\subset \supp{\mathcal V}\right\}$$
\item For $V_i\in\mathcal U$ compute a cover $\mathcal V_i\subset \mathcal Y$ of $f(V_i)$. This step is often realized by means of the interval arithmetic \cite{M}.
In most cases it is very difficult to compute $\mathcal V_i$ as a minimal cover of $f(V_i)$.
\item Define a cover of $f(U)$ by $\bigcup \mathcal V_i$.
\end{itemize}
The above considerations lead to a very natural definition of the combinatorial representation of a map. 
\begin{definition}
Let $f\colon X\to Y$ be a map and let us fix covers $\mathcal X,
\mathcal Y$ of $X$ and $Y$, respectively.

We will say that a combinatorial map $\mathcal F\colon\mathcal X\multimap\mathcal Y$ is a combinatorial representation of $f$ if for $x\in X$ and for all $U\in\mathcal X$ such that $x\in U$ holds $f(x)\in \supp{\mathcal F(U)}$.
\end{definition}
For a map $f\colon X\to X$ and a cover $\mathcal X$ of $X$ it is natural to encode a combinatorial representation of $f$ as a directed graph. A finite directed graph is a pair $\mathcal G = \graph{V}{E}$, where
$\mathcal V$ is a finite set whose elements are called vertexes, and $\mathcal E \subset \mathcal V\times \mathcal V$ is a set of selected
(ordered) pairs of vertexes. The elements of $\mathcal E$ are called edges. 
\begin{definition}
 We say that a graph $\mathcal G=\graph{V}{E}$ is a representation of a combinatorial map $\mathcal F\colon \mathcal U\multimap\mathcal U$ if
  $\mathcal V=\mathcal U$ and
 \begin{equation*}
  \mathcal E =  \{(U,V)\in\mathcal U\times \mathcal U : V\in \mathcal F(U)\}.
 \end{equation*}
\end{definition}

We will need a notion of the transposed graph and outgoing edges. For a graph $\mathcal G=\graph{V}{E}$ by $\mathcal G^T=\graph{V'}{E'}$ we will denote a graph in which $\mathcal V'=\mathcal V$ and
$$\mathcal E'=\left\{\right(U,W):(W,U)\in\mathcal E\}.$$

For a vertex $V$ in a graph $\mathcal G=\graph{V}{E}$ by $\mathrm{out}(V,\mathcal G)$ we will denote the set of vertexes $$\mathrm{out}(V,\mathcal G) = \{P\in\mathcal V : (V,P)\in\mathcal E\}.$$

\subsection{Enclosure of an attractor.}
In this section we will describe an algorithm that we used to enclose an attractor. There are already existing packages for rigorous enclosures of invariant objects, see for example \cite{GAIO}.
We used, however, our own implementation dedicated for this purpose which seems to be easier and makes the software independent on the other libraries. We will use the fact that what we want to enclose is an attractor.

The idea is very simple. First, we fix a cover $\mathcal X$ of the observed attracting domain. Then we choose a box $V$ from observed attracting region, and we enclose its forward trajectory using the sets from the cover $\mathcal X$ as long as the trajectory does not leave $\supp{\mathcal X}$. Since the cover $\mathcal X$ is finite by its definition, after a finite number of steps there are no new sets in the cover of the trajectory of $V$ and we can stop the procedure. This approach generates an enclosure of some invariant set, that is expected to be an attractor.

The Algorithm~\ref{alg:enclosure} summarizes the above considerations.

\begin{algorithm}[htbp]\linesnumbered\label{alg:enclosure}
 \caption{enclose forward trajectory}
 \KwIn{$V$ : element of $\mathcal X$,\newline
         $\mathcal G=\graph{V}{E}$ : graph,\newline
         $f$ : a map
         }
  \KwResult{$\mathcal G$ - graph}
  \KwData{
    $\mathcal U, \mathcal F$ : subsets of $\mathcal X$,\newline
    $\mathcal G'=(\mathcal V',\mathcal E')$ : graph
  }
  \Begin{
  $(\mathcal V,\mathcal E)\leftarrow(\emptyset,\emptyset)$\;
  $\mathcal U\leftarrow\{V\}$\;

  \Repeat{$\mathcal U=\emptyset$}
  {
    $(\mathcal V',\mathcal E')\leftarrow(\emptyset,\emptyset)$\;
   \ForEach{$U\in \mathcal U$}{
        \If{$f(U)\not\subset\supp{\mathcal X}$}{\Return \textbf{Failure};}
        $\mathcal F\leftarrow [f(U)]_{\mathcal X}$; // computed cover of $f(U)$\;
        $\mathcal V'\leftarrow \mathcal V'\cup\{U\}\cup \mathcal F$\;
        \ForEach{$F\in \mathcal F$}{
          $\mathcal E' \leftarrow\mathcal E'\cup\{(U,F)\}$\;
        }
   }
   $\mathcal U\leftarrow \mathcal V'\setminus(\mathcal V\cup\mathcal U)$\;
   $(\mathcal V,\mathcal E)\leftarrow(\mathcal V,\mathcal E)\cup(\mathcal V',\mathcal E')$\;
  }
  \Return{$\mathcal G$};
}
\end{algorithm}

\begin{lemma}
Fix a cover $\mathcal X$ of $X$. Assume Algorithm~\ref{alg:enclosure} is called with its arguments $V$, $\mathcal G=\graph{V}{E}$ and $f$.
If the algorithms stops without the \textbf{Failure} result then
\begin{enumerate}
  \item[(i)] the combinatorial map $\mathcal F\colon\mathcal V\multimap \mathcal V$ encoded by the graph $\mathcal G=\graph{V}{E}$ is combinatorial representation of $f|_{\supp{\mathcal V}}$, with the natural cover $\mathcal V$ of $\supp{\mathcal V}$,
 \item[(ii)] $\supp{\mathcal V}$ is a positive invariant set for $f$, i.e. $f(\supp{\mathcal V})\subset \supp{\mathcal V}$,
\item[(iii)] $\bigcup_{i=0}^\infty f^i(V)\subset \supp{\mathcal V}$.
\end{enumerate}
\end{lemma}
\begin{proof}
First observe that the algorithm always stops. By contradiction, assume it does not stop. This means that in each iteration of the main repeat-until loop (lines 4-15) the set $\mathcal U$ is not empty. This implies, that the number of elements in the set $\mathcal V$, which is updated in line 14, is strictly increasing. But this is a contradiction with $\mathcal V\subset \mathcal X$, which is a finite set by the definition of cover.

Let $\mathcal G$ be the result returned by the algorithm and let $\mathcal F\colon\mathcal V\multimap\mathcal V$ be the combinatorial map represented by the graph $\mathcal G$.

We are proving the assertion (i). Let $x\in\supp{\mathcal V}$ and let $V\in \mathcal V$ be such that $x\in V$. Observe, that if $\mathcal F(V)$ is not empty then $f(V)\subset\supp{\mathcal F(V)}$ - as guaranteed in lines 9, 11, 12 and 14 of the algorithm. Hence, in this case the assertion follows. To finish this step we observe, that the main repeat-until loop stops when the set $\mathcal U$ is empty. But this set contains of exactly these sets $U\in\mathcal V$ (lines 10 and 13) for which the value $\mathcal F(U)$ is not computed yet. Hence, after the algorithm stops, for each $V\in\mathcal V$, $\supp{\mathrm{out}(V,\mathcal G)}$ contains a computed enclosure of $f(V)$.

Assertions (ii) follows from (i) and the definition of combinatorial representation.

To prove (iii) it is enough to observe that $V\subset\supp{\mathcal V}$ as guaranteed in line 3, 10 and 14 of the algorithm. Then the assertion follows from (ii).
\end{proof}

Using Algorithm~\ref{alg:enclosure} one can compute combinatorial representation of a map restricted to some positive invariant set, which is expected to be an attractor. For further applications this graph should be as small as possible and it should still contain enclosure of the invariant set of $f$ restricted to $\supp{\mathcal V}$. Therefore it is important to choose the starting box $V$ such that it contains points from the attractor. An easy way to do that is to take a point from the observed basin of attraction and compute non-rigorously its forward trajectory for some long time. Then we can choose a set $V$ from the cover that contains the last computed point.

\subsection{Computing of coordinate systems.}
In this section we assume that the graph $\mathcal G=\graph{V}{E}$ is a combinatorial representation of $f$ resulting from Algorithm~\ref{alg:enclosure}. We present a heuristic method for computing of coordinate systems $C_V$ at the vertex $V$,  for $V\in \mathcal V$. We assume that at the beginning the coordinate systems are not set, which we will denote by $C_V$=NULL.

The method consists of the following steps.

\vskip\baselineskip
\noindent\textbf{Step 1. Finding cycles in the graph.} First we try to find (non-rigorously) as much as possible of periodic points of $f$ in $\supp{\mathcal V}$. This is important, because given a period $p$ point, say $x$, a very natural choice of a coordinate system in a vertex containing this point is the Jordan basis of $Df^p(x)$. Of course it may happen that several periodic points belong to the same vertex. In this case we are choosing one from these points of the lowest principal period.

More formally; let us fix a positive integer $N$. This is an upper bound for the highest period of orbits we will search for. Since $\mathcal G$ is a combinatorial representation of $f$, all periodic points of $f$ must belong to some cycles in the graph $\mathcal G$. Let $\mathcal P_i$, $i=1,\ldots N$ be a set of vertexes $V$ such that there exists an $i$-periodic cycle in the graph $\mathcal G$ through $V$ and no lower period cycle containing $V$ exists. The sets $\mathcal P_i$ can be computed by means of standard graph algorithms and we omit details here.

\vskip\baselineskip
\noindent\textbf{Step 2. Refining cycles to periodic points.}
In this step we compute the sets of points $P_i$, $i=1,\ldots,N$ that are good approximations of periodic points.
For $V\in\mathcal P_i$ we proceed as follows. Let $x$ be an approximate center of $V$. The point $x$ is a seed point for the standard Newton method for zero finding problem applied to the map $g(x)=f^i(x)-x$. If the Newton method converges to a point $y\in \supp{\mathcal V}$ which is of principal period $i$ then we insert the point $y$ to the set $P_i$.

We would like to comment here, that we do not try to find all periodic orbits of low periods for the map $f$. Using our approach one can detect most of them if the set $\mathcal V$ consists of small enough sets. Clearly, a larger number of detected periodic points is better for setting the coordinate systems but the task of finding all low period orbits is not a goal for us.

\vskip\baselineskip
\noindent\textbf{Step 3. Setting coordinate systems in a neighborhood of periodic points.}
The goal of this step is to assign a coordinate system to vertexes that contain computed approximate periodic point.
For $y\in P_i$ we proceed as follows. Let $V\in\mathcal V$ be a vertex such that $y\in V$. Such a vertex exists since we know that $y\in\supp{\mathcal V}$ -- as guaranteed by the previous step.
\begin{itemize}
\item If the vertex $V$ has already assigned a coordinate system then we proceed with the next point from $P_i$.
\item We compute an approximate derivative $Df^i(y)$ and an approximate matrix $M_y$ of normalized  eigenvectors. The columns of $M_y$ are sorted by decreasing absolute value of the associated eigenvalue.
\item We assign to the vertex $V$ a coordinate system $V_C = M_y$.
\end{itemize}
We would like to comment here that the sets $P_i$ are proceeded by increasing index $i$. This guarantees that a vertex $V$ will have assigned coordinate system from some periodic point of the lowest principal period. Note, the resulting coordinates depend on the order chosen in the set $P_i$.

\vskip\baselineskip
\noindent\textbf{Step 4. Spreading coordinates from periodic points.}
The main idea is to propagate the coordinate systems from periodic points by the action of the derivative. It is well known that the direct propagation of coordinates by $Df$ in floating point arithmetic usually results in collapsing of these coordinates to a singular matrix. Therefore we will involve some orthonormalization process.

The Algorithm~\ref{alg:spreadCoordinates} describes the procedure.
\begin{algorithm}[htbp]\linesnumbered\label{alg:spreadCoordinates}
 \caption{spread coordinates from periodic points}
  \KwIn{$\mathcal G=\graph{V}{E}$ : graph,\newline
         $\mathcal S$ : subset of $\mathcal V$,\newline
         $k>1$ : integer
         }
  \KwResult{$\{C_V\}_{V\in\mathcal V}$ - coordinate systems at each vertex in the graph $\mathcal G$}
  \KwData{
  $u$ : vector;\newline
  $M$ : stack of float matrices;\newline
  $C$ : matrix;\newline
  $\mathcal P$ : subset of $\mathcal V$;\newline
  $i$ : integer
  }
  \Begin{
  \Repeat{$S=\emptyset$}
  {
    $\mathcal P\leftarrow \emptyset$\;
   \ForEach{$V\in \mathcal S$}{
        $u\leftarrow \mathrm{centre}(V)$\; // we will save on the stack derivatives along the trajectory\;
        $C\leftarrow C_V$\;
        \For{$i=1$ to $k$}{
          push($Df(u)$,$M$)\;
          $C\leftarrow Df(u)\cdot C$\;
          $u\leftarrow f(u)$;\
        }
        // then we compute coordinate system at the image\;
        $C\leftarrow \mathrm{GrammSchmidtOrthonormalization}(C)$\;
        \For{$i=1$ to $k-1$}{
          $C\leftarrow \left(\mathrm{top}(M)\right)^{-1}\cdot C$\;
          pop($M$)\;
        }
        pop($M$); // to free the memory\;
        $C\leftarrow\mathrm{normalizeColumns}(C)$\;
        \ForEach{$W\in\mathrm{out}(V,\mathcal G)$}
        {
          \If{$C_W\neq$NULL}{
            $C_W\leftarrow C^{-1}$\;
             $\mathcal P\leftarrow\mathcal P\cup\{W\}$\;
          }
        }
    }
    $\mathcal S\leftarrow\mathcal P$\;
  }
  }
\end{algorithm}

\begin{lemma}\label{lem:spreadCoordinates}
 Assume the Algorithm~\ref{alg:spreadCoordinates} is called with its arguments $\mathcal G=\graph{V}{E}$, $\mathcal S\neq\emptyset$ and $k>1$. Assume $C_V\neq NULL$ for $V\in \mathcal S$. If $\mathcal G$ consists of only one strongly connected component then the algorithm stops and $C_V\neq NULL$ for all $V\in\mathcal V$.
\end{lemma}
\begin{proof}
  We will prove that the algorithm always stops. Assume it is not the case. Then in each iteration of the main loop the set $\mathcal P$ is non-empty. But this set contains (lines 3 and 19-22) exactly these vertexes for which the coordinate system is assigned in the current iteration of the loop. Hence, the set of vertexes at which we have assigned the coordinate system will increase in each iteration. This contradicts the fact, that $\mathcal V$ is a finite set.

  We will prove that the algorithm sets coordinates in each vertex. Assume it is not the case.
  Take $\mathcal G^T=(\mathcal V',\mathcal E')$. Let $V$ be a vertex such that $C_V$=NULL. Then for $U\in\mathrm{out}(V,\mathcal V')$ holds $C_U=NULL$. Otherwise, $U$ will belong to the set $\mathcal S$ in some iteration and then the coordinates will be set in $V$. Denote $\mathcal V_1=\mathrm{out}(V,\mathcal G^T)$.

   In a recurrent way we can define $$\mathcal V_k = \bigcup\left\{\mathrm{out}(U,\mathcal G^T) : U\in \mathcal V_{k-1}\right\}.$$ The same reasoning proves that $C_U$=NULL for $U\in\mathcal V_k$.

  Since $\mathcal G$ consists of only one strongly connected component, $\mathcal G^T$ does too. Hence, for some $K>0$ we have $V_K=\mathcal V$ which contradicts the assumption that $C_U\neq$NULL for $U\in\mathcal S$ and $\mathcal S\neq\emptyset$.
\end{proof}

The input argument $\mathcal S$ is the set of vertexes that contain periodic points and for which the coordinate systems were already computed in Step 3. As we have seen, lines 5-18 are not important for the correctness of the Algorithm~\ref{alg:spreadCoordinates}. In these lines a heuristic method for computing of coordinate systems is proposed. We propagate the actual coordinate system $C_V$ by the derivative along a short part ($k$ iterations) of the trajectory of the center of $V$. Then we perform orthonormalization of columns and propagate the obtained matrix backwards $k-1$ times. This heuristic can work because for the inverse map the stable directions are attracting for $Tf$ and after few steps we get some approximation of stable directions at $f(\mathrm{centre}(V))$.

In this step we also assumed that the graph $\mathcal G$ consists of only one strongly connected component. It is not a restriction, since the method is dedicated for attractors. In practice, Algorithm~\ref{alg:enclosure} always returns such a graph if the initial set $V$ is chosen carefully.  Formally, it is easy to verify that a graph consists of one strongly connected component using for example Tarjan's algorithm \cite{Ta}.

\begin{remark}
The result of the algorithms described in steps 2,3,4 depend on the order chosen to index the sets of vertexes and points ($\mathcal P_i$, $P_i$  and $\mathcal S$).
Moreover, these algorithms are not deterministic, when implemented in parallel mode. The coordinate system at each vertex is computed by the first thread that reaches this vertex.

Anyway, the conclusion of Lemma~\ref{lem:spreadCoordinates} holds; we have computed some coordinate system at each vertex and we can start verification of the cone conditions.
\end{remark}

\subsection{Verification of the cone conditions.}
This step is quite straightforward. Let $\mathcal G=\graph{V}{E}$ be a graph which encodes a combinatorial representation of a map $f$. Assume we have computed some coordinate systems $\{C_V\}_{V\in\mathcal V}$. The goal is to check whether the assumptions of Theorem~\ref{thm:verifyConesByQuadForm} hold true for $f$. This can be done by simply verifying the positive definiteness of some interval matrices corresponding to every edge in $\mathcal E$.

\begin{algorithm}[htbp]\linesnumbered\label{alg:verifyConeConditions}
 \caption{verify cone conditions}
 \KwIn{$\mathcal G=\graph{V}{E}$ : graph;\newline
         $Q$ : matrix;\newline
         $f$ : map;
         }
   \KwResult{$\mathcal U$ : subset of $\mathcal V$;}
  \KwData{
    $V,W$ : element of $\mathcal V$;\newline
    $D, M, A$ : interval matrix;\newline
  }
  \Begin{
   $\mathcal U\leftarrow \emptyset$\;
   \ForEach{$V\in \mathcal V$}{
      $D\leftarrow [Df(V)]$; // computed interval enclosure of $Df(V)$\;
      \ForEach{$W\in \mathrm{out}(V,\mathcal G)$}{
        $M\leftarrow [C_W D C_V^{-1}]$\;
        $A\leftarrow M^T\cdot Q\cdot M - Q$\;
        \If{cannot verify that $A$ is positive definite}{
          $\mathcal U = \mathcal U\cup\{V\}$\;
        }
    }
  }
  \Return{$\mathcal U$;}
}
\end{algorithm}

We have the following obvious lemma.
\begin{lemma}\label{lem:verifyConeConditions}
 Assume the Algorithm~\ref{alg:verifyConeConditions} is called with its arguments $\mathcal G=\graph{V}{E}$, $Q$ and $f$. Assume that all the matrices $\{C_V\}_{V\in\mathcal V}$ are invertible. The algorithm always stops and returns a set $\mathcal U$ such that for $V\in\mathcal V\setminus\mathcal U$ and $W\in\mathrm{out}(V,\mathcal G)$ the interval matrix
\begin{equation*}
  [C_WDf(V)C_V^{-1}]^TQ[C_WDf(V)C_V^{-1}]-Q 
\end{equation*}
is positive definite.
\end{lemma}

A direct consequence of the above Lemma and Theorem~\ref{thm:verifyConesByQuadForm} is the following
\begin{corollary}
 The same assumption as in Lemma~\ref{lem:verifyConeConditions}. Assume $\mathcal G$ encodes a combinatorial representation of the map $f$ and let $Q$ be quadratic form as in (\ref{eq:Qform}). If the Algorithm~\ref{alg:verifyConeConditions} called with its arguments $\mathcal G=\graph{V}{E}$, $Q$ and $f$ returns an empty set $\mathcal U$ then $f$ is uniformly hyperbolic on $\mathrm{Inv}(f,\supp{\mathcal V})$.
\end{corollary}

\section{Applications.}\label{sec:applications}
In this section we will show how the method introduced in the last section works in a very easy example. Then we will give a proof of Theorem~\ref{thm:main}.

\subsection{Toy example - Smale map.}
\begin{figure}[htbp]
 \centerline{
    \includegraphics[width=0.5\textwidth]{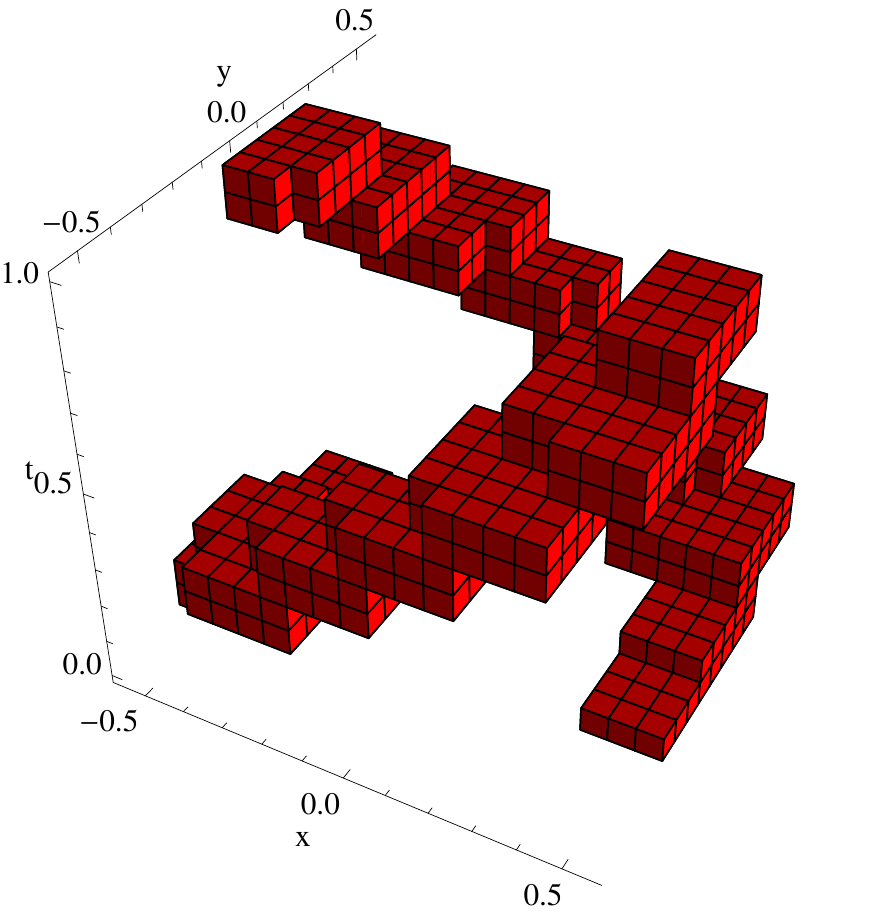}
    \includegraphics[width=0.5\textwidth]{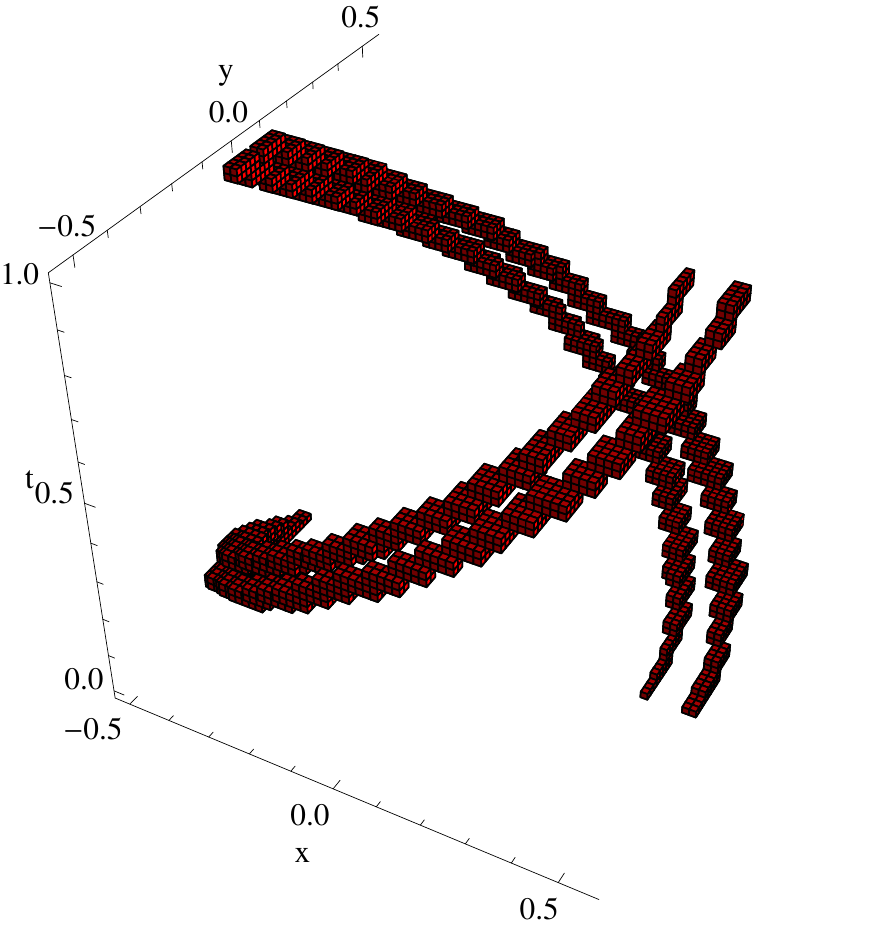}
  }
\caption{Enclosure of the Smale attractor for the resolutions $k=4$ and $k=6$. Note, the planes $t=0$ and $t=1$ should be identified.\label{fig:smale}}
\end{figure}

A natural choice to test the algorithms and implementation is a very simple example. Here we have chosen the well known Smale map $s\colon\mathbb R^2\times [0,1]\to\mathbb R^2\times [0,1]$
\begin{equation}
 s(x,y,t)=(0.1x+0.5\cos(2\pi t),0.1y+0.5\sin(2\pi t),2t\ \mathrm{mod}\ 1).
\end{equation}
The map is proved to be uniformly hyperbolic just by simple hand-made calculations. It is easy to see that the set $[-1,1]\times[-1,1]\times [0,1]$ is positive invariant for $s$. Define a cover of this attracting domain of the form
\begin{equation*}
 \mathcal X_k=\left\{
  \frac{
    \left[a,a+1\right]\times
    \left[b,b+1\right]\times
    \left[c,c+1\right]
  }{2^k}
  : a,b,c =-2^k,\ldots,2^k-1, c\geq0
  \right\}.
\end{equation*}
We will call the number $k$ in $\mathcal X_k$ the \emph{resolution}.
After application of Algorithm~\ref{alg:enclosure} we got an enclosure of the attractor consisting of $567$ boxes for $k=4$ and $3333$ boxes for $k=6$. These enclosures are shown in Fig.~\ref{fig:smale}.

We completed the proof of uniform hyperbolicity for the resolution $k=4$. We run the algorithm for finding cycles in the graph and refining them to periodic points with the maximum period set to $N=3$. They returned $1$ candidate for fixed point, a period $2$ orbit and two period $3$ orbits.

Then we applied the algorithm for spreading coordinate systems with $k=2$ (this is the number of forward iterates when propagate the coordinate systems) and verification of the cone conditions. We run Algorithm~\ref{alg:verifyConeConditions} with the quadratic form
\begin{equation*}
 Q=\begin{bmatrix}
    1 & 0 & 0\\
    0 & -1 & 0\\
    0 & 0 & -1
   \end{bmatrix}
\end{equation*}
and it returned an empty set $\mathcal U$ of unverified vertexes. The program executes within less than 1 second on a laptop-type computer.

\subsection{Application to the system of coupled van der Pol oscillators -- proof of Theorem~\ref{thm:main}.}
Consider the system (\ref{eq:vpo}) with the parameter values (\ref{eq:parameters}). Let $P$ be the Poincar\'e map as defined in (\ref{eq:poincare}). Let us perform a linear change of coordinates
\begin{equation*}
 x_0 = x/r_0,\quad x_1=(u-c_{ux})/r_1,\quad x_2=y-c_{yx}x-c_{yu}u,\quad x_3=v-c_{vx}x-c_{vu}-c_{vy}y,
\end{equation*}
where
\begin{equation*}
  \begin{array}{ll}
  c_{ux}=0.438, & c_{yx} = -0.042,\\
  c_{yu}=0.226, & c_{vx} = -0.218,\\
  c_{vu} = 0.029, & c_{vy} = -0.118,\\
  r_0 = 0.812, & r_1 = 0.721.
  \end{array}
\end{equation*}
As it is observed in \cite{KSe}, in these coordinates the attractor is more aligned to the coordinate axes. In the sequel we will consider the Poincar\'e map $\tilde P = L^{-1} P L$, where $L$ is the matrix of coordinate change.

After few attempts we found that the resolution $k=14$ is enough for verification of the hyperbolicity of the observed attractor for $\tilde P$.  This means that we have used a uniform grid in some bounded domain of the form
\begin{equation*}
 \mathcal X=\left\{
  \frac{
    \left[a,a+1\right]\times
    \left[b,b+1\right]\times
    \left[c,c+1\right]\times
    \left[d,d+1\right]
  }{2^k}
  : a,b,c,d =-2^{k+1},\ldots,2^{k+1}-1
  \right\}.
\end{equation*}
In Table~\ref{tab:swData} we give the results of applying of the algorithms for generating enclosure, setting of coordinate systems and verification of strong hyperbolicity.
\begin{table}
\begin{center}
 \begin{tabular}{|c|c|c|}
  \hline
  algorithm & wall time (h:mm) & comments\\
  \hline
  enclosure of attractor, &  & \\
  Algorithm~\ref{alg:enclosure} & 2:16 on 224 CPUs & $7\,970\,392$ boxes\\
  \hline
  finding of cycles in graph,  & &\\
  max period $6$ & 0:58 on 32 CPUs & $2190$ cycles found\\
  \hline
  finding of periodic points,& &\\
  max period $6$ & 0:06 on 32 CPUs & $105$ points found\\
  \hline
  computing of coordinate systems, & & \\
  Step 4 and Algorithm~\ref{alg:spreadCoordinates}& 6:59 on 32 CPUs & parameter $k=2$ \\
  \hline
  verification of cone condition, & & empty set \\
  Algorithm~\ref{alg:verifyConeConditions}& 4:24 on 224 CPUs & of unverified boxes \\
  \hline
 \end{tabular}\vskip .3cm
\end{center}
\caption{Setting, results and time of computation of the algorithms applied to the map $\tilde P$. \label{tab:swData}}
\end{table}
We run Algorithm~\ref{alg:verifyConeConditions} with the quadratic form
\begin{equation*}
 Q=\begin{bmatrix}
    1 & 0 & 0 & 0 \\
    0 & -1 & 0 & 0 \\
    0 & 0 & -1 & 0 \\
    0 & 0 & 0 & -1
   \end{bmatrix}
\end{equation*}
and it returned an empty set of unverified vertexes. Hence, the uniform hyperbolicity of the set $\mathrm{Inv}(\tilde P,\supp{\mathcal V})$ is verified. 
To complete the proof of Theorem~\ref{thm:main} we used the CHOMP library \cite{CHOMP} in order to compute homology groups of the cubical set represented by $\mathcal V$. The program verified that $H_0(\mathcal V)=\mathbb Z$, hence the set $\supp{\mathcal V}$ is connected. The group $H_1(\mathcal V)$ also was nonzero, hence the set is non-contractible.

This finishes the proof of Theorem~\ref{thm:main}.

\subsection{The existence of a fixed point and a period two orbit -- proof of Theorem~\ref{thm:nontrivAttractor}.}
As a byproduct of the proof of Theorem~\ref{thm:main} we got a good numerical approximations of over one hundred periodic points. To prove that the set $\mathrm{Inv}(\tilde P,\supp{\mathcal V})$ is nontrivial we used the interval Newton operator \cite{M} as a tool for verification of the existence and local uniqueness of zeros of maps. We have the following lemma.
\begin{lemma}\label{lem:newton}
 Assume $f$ is a smooth map in a neighborhood of an interval set $X\subset\mathbb R^n$ and fix $\bar x\in\mathrm{int} X$. If the interval Newton operator
\begin{equation*}
 N(\bar x,X,f) = \bar x - [Df(X)]_I^{-1}(X-\bar x)\subset\mathrm{int}(X)
\end{equation*}
then the map $f$ has unique zero in the set $X$. Moreover, if $x_*\in X$ is the unique zero of $f$ in $X$ then $x_*\in N(\bar x,X,f)$.
\end{lemma}

\begin{proof}[Proof of Theorem~\ref{thm:nontrivAttractor}]
A good numerical approximation of a fixed point for $\tilde P$ is
\begin{equation*}
\bar x = \begin{bmatrix}
    -0.2217320903523841312\\
    0.97492087818262221051\\
    0.018539814418960732759\\
    -0.03127925657964333925
     \end{bmatrix}^T
\end{equation*}
-- this is a point returned from the step of finding periodic points in $\supp{\mathcal V}$. Let $X$ be a ball of radius $10^{-10}$ in the maximum norm centered at $\bar x$. We computed the interval Newton operator $N(\bar x, X,\tilde P-\mathrm{Id})$ and we got that
\begin{equation*}
 \|N(\bar x, X, \tilde P-\mathrm{Id}) - \bar x\|_0 \leq 2.5\cdot 10^{-12}.
\end{equation*}
Hence, $N(\bar x, X, \tilde P-\mathrm{Id})\subset \mathrm{int}(X)$ and the map $\tilde P-\mathrm{Id}$ has a unique zero in $X$. Clearly this is a fixed point for $\tilde P$. We also verified that $N(\bar x, X, \tilde P-\mathrm{Id})\subset \supp{\mathcal V}$ and this finishes the proof that $\tilde P$ has a fixed point $\supp{\mathcal V}$.

In a similar way we verified that the attractor contains a period two orbit. It is enough to show that the map $F(x,y) = (P(x)-y,P(y)-x)$ has a zero for some $(x_*,y_*)$ and $x_*\neq y_*$. We have a good candidate for periodic points
\begin{equation*}
 \bar x = \begin{bmatrix}
           1.0708168699357167158\\
           0.093670228872794565237\\
           0.015571609254819055212\\
          0.024616361886390349154
          \end{bmatrix}^T,\qquad
\bar y = \begin{bmatrix}
          -0.22099700889399603573\\
          -0.706209321733317795\\
          -0.044189820507618084004\\
          -0.006451355101601539937
         \end{bmatrix}^T.
\end{equation*}
Let $Z\subset \mathbb R^8$ be a ball of radius $10^{-10}$ in the maximum norm centered at $\bar z = (\bar x, \bar y)$. We verified that
\begin{equation*}
 \|N(\bar z, Z, F) -\bar z\|_0 \leq 2\cdot 10^{-12}.
\end{equation*}
This proves that $F$ has a unique zero $(x_*,y_*)\in Z$. Since $\|\bar x - \bar y\|_0\geq 2\cdot10^{-10}$ we have $x_*\neq y_*$ and both $x_*, y_*$ are period two points for $\tilde P$. We also verified that $(x_*,y_*)\in N(\bar z, Z, F)\subset \supp{\mathcal V}\times \supp{\mathcal V}$ and this completes the proof of Theorem~\ref{thm:nontrivAttractor}.
\end{proof}

\section{Conclusions.}
In this paper we proposed a method for rigorous verification of uniform hyperbolicity for maps. The method has been successfully applied to a Poincar\'e map of a non-autonomous ODE in four dimensions. We believe the method can be applied to other systems.

\begin{figure}[htbp]
 \centerline{
    \includegraphics[width=0.75\textwidth]{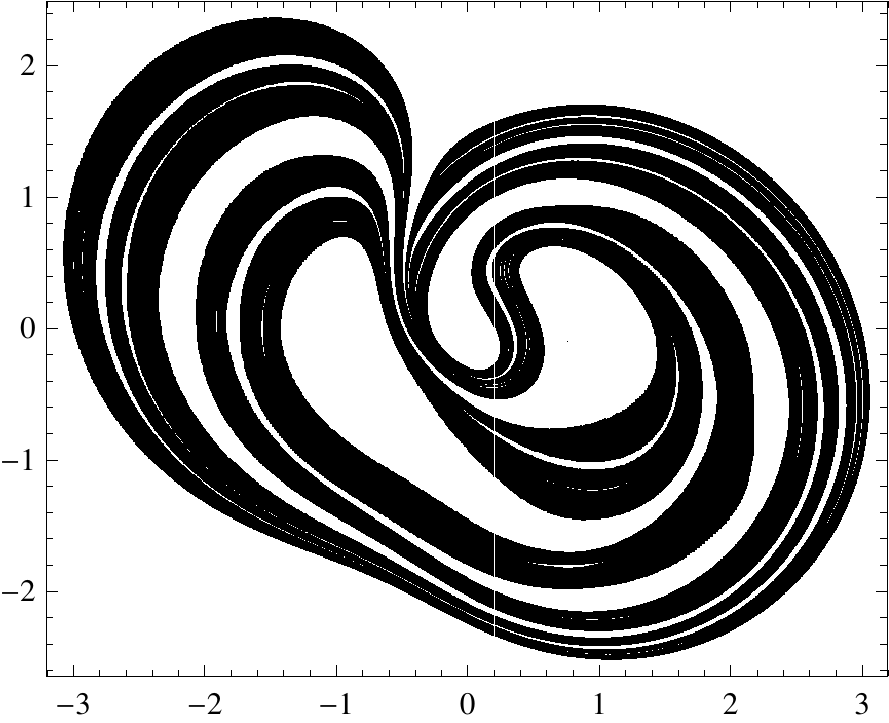}
  }
\caption{Enclosure of the attractor of the Poincar\'e map for the system (\ref{eq:plykin}) with parameter values $K=1.9$, $\varepsilon=0.72$ and with the resolution $k=9$. \label{fig:plykin}}
\end{figure}

Kuznetsov \cite{K2} proposed also a non-autonomous system on the plane which apparently possesses an attractor of Plykin type in the Poincar\'e map. The equations are quite complicated and given by
\begin{equation}\label{eq:plykin}
  \begin{array}{rcl}
 \dot x &=& -2\varepsilon y^2\Omega_1(x,y,t)
      \left(
          \cos\left(\frac{\pi}{4}\cos\frac{\pi}{2}t\right)
          -x\sin\left(\frac{\pi}{4}\cos\frac{\pi}{2}t\right)
      \right) + \\
      &&  Ky\Omega_2(x,y,t)\left(
          \cos\left(\frac{\pi}{4}\sin\frac{\pi}{2}t\right)
          -x\sin\left(\frac{\pi}{4}\sin\frac{\pi}{2}t\right)
      \right)\sin\frac{\pi}{2}t,\\\\
  \dot y &=& \varepsilon y\Omega_1(x,y,t)
      \left(
        2x\cos\left(\frac{\pi}{4}\cos\frac{\pi}{2}t\right)+(1-x^2+y^2)\sin\left(\frac{\pi}{4}\cos\frac{\pi}{2}t\right)
      \right) -\\
    && K\Omega_2(x,y,t)\left(
      x\cos\left(\frac{\pi}{4}\sin\frac{\pi}{2}t\right)+\frac{1}{2}(1-x^2+y^2)\sin\left(\frac{\pi}{4}\sin\frac{\pi}{2}t\right)
    \right)\sin\frac{\pi}{2}t,
  \end{array}
\end{equation}
where
\begin{eqnarray*}
 \Omega_1(x,y,t) &=& \frac{2x\cos\left(\frac{\pi}{4}\cos\frac{\pi}{2}t\right)+(1-x^2-y^2)\sin\left(\frac{\pi}{4}\cos\frac{\pi}{2}t\right)}{\left(1+x^2+y^2\right)^2},\\
 \Omega_2(x,y,t) &=& \frac{-2x\sin\left(\frac{\pi}{4}\sin\frac{\pi}{2}t\right)+(1-x^2-y^2)\cos\left(\frac{\pi}{4}\sin\frac{\pi}{2}t\right)}{\left(1+x^2+y^2\right)^2} + \frac{1}{\sqrt 2}.
\end{eqnarray*}
Kuznetsov \cite{K2} did a numerical study of this system and observed that for a wide range of parameter values around $K=1.9$ and $\varepsilon=0.72$ the Poincar\'e map has a hyperbolic attractor.
A rigorous enclosure of the attractor for $K=1.9$, $\varepsilon=0.72$ resulting from Algorithm~\ref{alg:enclosure} with the resolution $k=9$ is shown in Fig.~\ref{fig:plykin}. The resolution $k=9$ was not enough to perform the step of verification of uniform hyperbolicity of this map. We also enclosed this attractor with resolution $k=12$ but again this was not enough to verify the hyperbolicity. Further subdivision was not possible due to memory limitations in the computer that we used.

Computations for the system (\ref{eq:plykin}) are much harder than for (\ref{eq:vpo}). First - the vector field is a quite large expression and even if the system is lower dimensional, it is more complicated for rigorous integration. Second, and most important; as easily seen in Fig.~\ref{fig:plykin}, the attractor is very thick and harder to cover than the almost one dimensional Smale solenoid. An approximate Hausdorff dimension reported in \cite{K2} is ca $1.7$. We reached the limit of available memory on our computer when enclosing the attractor with a thinner cover. We believe, however, that some further optimizations in data structures will help to prove hyperbolicity of this system. A
straightforward idea is to use non-uniform covers.

\subsection{General invariant sets.}

In this paper we verified hyperbolicity of attractors. But the method can be easily extended to invariant sets. The only part which requires some modification is the method for generating an enclosure of an invariant set. If the invariant set is an attractor, the strategy of ``inner'' enclosure is better - Algorithm~\ref{alg:enclosure}. For a general invariant set an ``outer'' approximation is necessary. We also implemented ``outer'' approximation for enclosing of general invariant sets but we do not report details here, since it is rather standard procedure. As a test case we proved that the well known H\'enon map \cite{H} 
\begin{equation*}
\mathcal H_{a,b}(x,y) = (1 + y- ax^2,bx)
\end{equation*}
is uniformly hyperbolic on the invariant set for the parameter values $a=5.4$ and $b=-1$. 

\begin{figure}[htbp]
 \centerline{
    \includegraphics[width=0.75\textwidth]{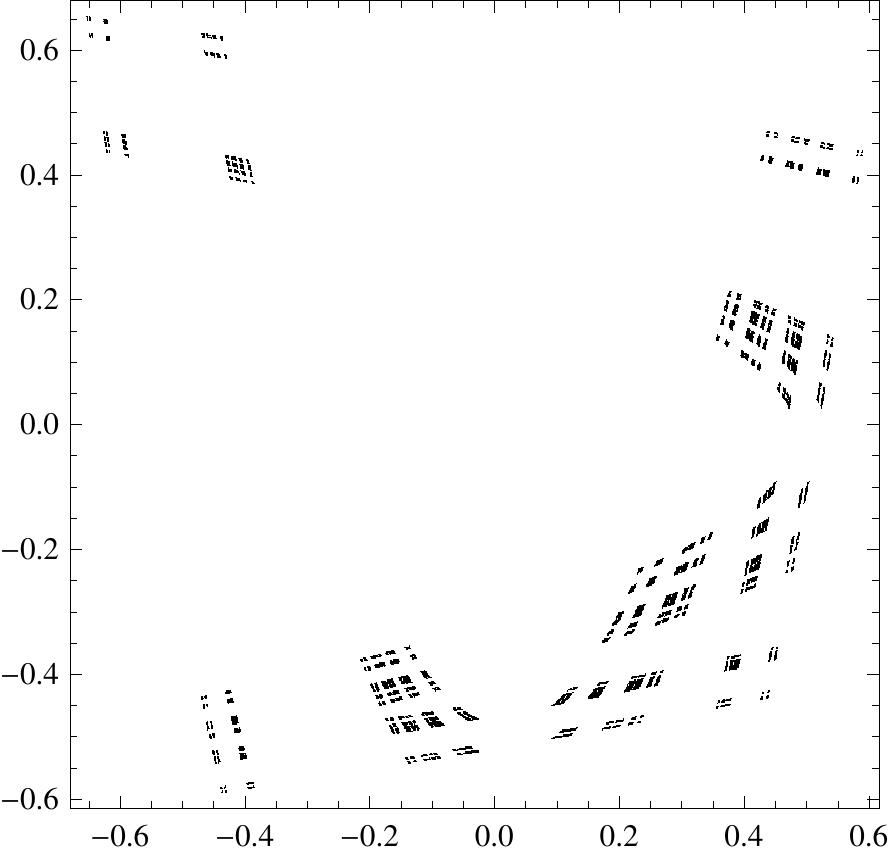}
  }
\caption{An enclosure of the invariant set for $\mathcal H_{a,b}$ map with the parameter values $a=5.4$, $b=-1$. \label{fig:henon}}
\end{figure}

As mentioned above, we used an ``outer'' strategy to enclose the invariant set -- the result is shown in Fig.~\ref{fig:henon} and it consists of $8832$ boxes. Then we applied algorithms for finding cycles and periodic points, computing of coordinate systems and verification of the cone condition with the quadratic form $$Q=\begin{bmatrix}1 & 0 \\ 0 & -1\end{bmatrix}.$$

The program executes within less than $1$ second on a laptop-type computer with the Intel Core 2 Duo 2GHz processor. A computer assisted proof for the same parameter values is presented in \cite{A,MT}. In \cite{MT} the authors report the time of computations less than $10$ seconds on a comparable CPU.

\subsection{Implementation.}
The C++ program for verifying the hyperbolicity has been implemented by the author and the source code is available at \cite{W}. The code is highly parallelized using the Open MP library supported by the compilers gcc-4.2 or newer. The code is written in a very generic way; it heavily uses template techniques. In fact, an input to the algorithms is an abstract map (template parameter) for which we assume that we know how to compute values and derivatives. In the case of Poincar\'e maps we used for this purpose the $C^0$ and $C^1$ solvers from the \cite{CAPD} library -- see also \cite{Z1} for an efficient algorithm for integration of first order variational equations.

All computations regarding Theorem~\ref{thm:main} and Theorem~\ref{thm:nontrivAttractor} were performed on a cluster with $7$ computers, each having $8$ Quad-Core AMD Opteron(tm)
8354, 2.2GHz processors with 64GB of RAM. During the computations the maximal usage of the memory was 32\%.


\begin{thebibliography}{CAPD}
\bibitem[A]{A} Z. Arai, \emph{On Hyperbolic Plateaus of the H\'enon map}, Experimental Mathematics, 16:2 (2007), 181--188.
%
\bibitem[CAPD]{CAPD} CAPD --- Computer assisted proofs in dynamics, a package
for
rigorous numerics, \verb| http://capd.ii.uj.edu.pl.|
%

\bibitem[CHOMP]{CHOMP} CHOMP -- Computational Homology Project,
{\tt http://chomp.rutgers.edu}.

\bibitem[GAIO]{GAIO} M. Dellnitz and O. Junge, The web page of GAIO project,
{\tt http://math-www.uni-paderborn.de/~agdellnitz/gaio}.

\bibitem[GH]{GH} J. Guckenheimer and P. Holmes, \emph{Nonlinear Oscillations,
 Dynamical Systems and Bifurcations of vector Fields}, Springer-Verlag,
    New York, Vol. 43 of Applied Math. Sciences (1983).

\bibitem[H]{H} M. H\'enon, \emph{A two-dimensional mapping with a strange
attractor}, Comm. Math. Phys. \textbf{50}  (1976), 69--77.
%
\bibitem[Hr1]{Hr1} S.L. Hruska, \emph{A Numerical Method for Constructing the Hyperbolic
Structure of Complex H\'enon Mappings}, Found. Comp. Math. 6, No. 4, (2006), 427--455.

\bibitem[Hr2]{Hr2} S.L. Hruska, \emph{Rigorous numerical models for the dynamics of complex H\'enon mappings on their
chain recurrent sets}, Discrete Contin. Dynam. Syst., 15 (2) (2006), 529--558.

\bibitem[K]{K} S.P. Kuznetsov, \emph{Example of a Physical System with a Hyperbolic Attractor of the Smale-Williams Type}, Phys. Rev. Lett., 95, 2005, 144101.

\bibitem[K2]{K2} S.P. Kuznetsov, \emph{A non-autonomous flow system with Plykin type attractor}, Communications in Nonlinear Science and Numerical Simulation, 14, 2009, 3487--3491,
 
\bibitem[KH]{KH} A. Katok and B. Hasselblatt, {\em Introduction to the Modern Theory of
Dynamical Systems}, Cambridge: Cambridge University Press; 1995.

\bibitem[KS]{KS} S.P. Kuznetsov and I.R. Sataev, {\em Hyperbolic attractor in a system of coupled non-autonomous van der Pol oscillators: Numerical test for expanding and contracting cones}, Phys. Lett. A 365, 97--104, (2007).

\bibitem[KSe]{KSe} S.P.Kuznetsov, E.P.Seleznev, \emph{A strange attractor of the Smale--Williams type in the chaotic dynamics of a physical system}, JETP 102, 2006, No. 2, 355--364.

\bibitem[KWZ]{KWZ} H. Kokubu, D. Wilczak and  P. Zgliczy\'nski,
    \emph{Rigorous verification of cocoon bifurcations in the Michelson system},
    Nonlinearity, \textbf{20} (2007), 2147--2174.

\bibitem[MT]{MT} M. Mazur, J. Tabor, \emph{Computational hyperbolicity}, preprint.

\bibitem[MTK]{MTK} M. Mazur, J. Tabor, P. Ko\'scielniak, \emph{Semi-hyperbolicity and hyperbolicity}, Disc. Cont. Dyn. Sys. 20, No. 4 (2008), 1029--1038.

\bibitem[M]{M} R.E. Moore, {\em Interval Analysis.} Prentice
Hall, Englewood Cliffs, N.J., 1966.

\bibitem[P]{P} R.V. Plykin, {\em Sources and sinks of A-diffeomorphisms of
surfaces,} Math. USSR Sb. 23(2):233-253 (1974).

\bibitem[PT]{PT} J. Palis and F. Takens, \emph{Hyperbolicity \& sensitive
chaotic dynamics at homoclinic bifurcations}, Cambridge studies in advanced
mathematics, vol. 35, Cambridge University Press, 1993

\bibitem[T]{T} W. Tucker, \emph{A Rigorous ODE solver and Smale's 14th Problem},
Found. Comput. Math., 2:1, 53--117, 2002.

\bibitem[Ta]{Ta} R. Tarjan, \emph{Depth-first search and linear graph algorithms}, SIAM Journal on Computing, Vol. 1 (1972), 146--160.

\bibitem[W]{W} D. Wilczak, {\tt http://www.ii.uj.edu.pl/\~{}wilczak},
a reference for auxiliary materials.

\bibitem[Z1]{Z1} P. Zgliczy\'nski,  \emph{$C^1$-Lohner algorithm}, Foundations
of Computational Mathematics, {\bf 2} (2002), 429--465.

\bibitem[Z]{Z} P. Zgliczy\'nski, \emph{Covering relations, cone conditions and
stable manifold theorem}, J. Diff. Eq., \textbf{246} (2009), 1774--1819.

\end{thebibliography}
\end{document}